\documentclass[a4paper,12pt]{article}
\usepackage[utf8]{inputenc}
\usepackage[english]{babel}
\usepackage{amssymb}
\usepackage{amsmath}
\usepackage{epsfig}
\usepackage{bm}
\usepackage{graphics}
\usepackage{pict2e}
\usepackage{xcolor}
\usepackage{tikz}

\usetikzlibrary{math}

\def\be{\begin{equation}}
\def\ee{\end{equation}}

\def\deg{{\rm deg}}

\newtheorem{thrm}{\bf Theorem}

\newtheorem{rmk}{\bf Remark}
\newtheorem{xmpl}{\bf Example}
\newtheorem{dfn}{\bf Definition}
\def\bl{\rule[-1mm]{2.4mm}{2.4mm}}

\def\RR{{\mathbb{R}}}
\def\CC{{\mathbb{C}}}
\def\ZZ{{\mathbb{Z}}}
\def\HH{{\mathbb{H}}}

\graphicspath{{figs/}{../Part I/figs/}}

\begin{document}

\title{Calculation of Stiefel fractions via Riemann thetas}
\author{Andrei Bogatyr\"ev and Timur Garaev
\thanks{Supported by INM RAS branch of MCFAM, agreement 075-15-2025-347.}
}
\date{}
\maketitle

\abstract{Stiefel fraction is a deformation of Zolotarev one which allows to solve certain three band extentions of Zolotarev approximation problem. We propose an efficient computational approach for Stiefel's solution based on genus two 
Riemann theta functions. This allows us to visualise the domain in the space of three band problems corresponding to  Stiefel's solutions.
}

\section*{Introduction}
Motivated by optimal synthesis of electrical filters Edward Stiefel has considered \cite{Sti61} a natural 
three band extension of the famous Zolotarev's third approximation problem \cite{Zol77}.

The set  $E:=E^-\sqcup E_1^+\sqcup E_2^+$ of active bands of a filter is a collection of three disjoint 
closed segments of the extended real axis $\hat{\RR}:={\mathbb R}\cup\infty$. Segment $E^-$ is called the stopband and $E^+:=E_1^+\sqcup E_2^+$ are two passbands. 

{\bf Multiband filter approximation problem}: Given the band set $E$ and positive integer $n$,  
minimize the \emph{goal function} among all real rational functions $R(x)$ of  bounded degree 
$\deg(R)\le n$ being the maximum of degrees of numerator and denominator for the 
uncancellable representation of the fraction $R(x)$.

Stiefel has considered the \emph{goal function 1}
\be
\lambda^2[R]:=\frac{\max_{x\in E^+}|R(x)|}{\min_{x\in E^-}|R(x)|}\to \min, 
\label{MW}
\ee
which is the three band extension of the third Zolotarev problem \cite{Zol77, Ach1, Ach2} and  
is typical for microwave filters. For the analogue filters they use another goal function which 
brings us to the extension of the forth Zolotarev problem.
Consider locally constant two-valued function $S_E(x)=\pm1$ when $x\in E^\pm$ (having the meaning of 
an ideal gain-frequency characteristic for the filter) and the related \emph{goal function 2}:
\be
\label{Analogue}
 \mu[R]:=||R-S_E||_{C(E)} := \max\limits_{x \in E}|R(x) - S_E(x)| \to \min. 
\ee

It can be shown \cite{B21} that both settings are equivalent. In particular, the optimal function $R_1$ for the \emph{goal function $1$} and normalization 
$\max_{x\in E^+}|R(x)|~\min_{x\in E^-}|R(x)|=1$ is related to optimizer  $R_2$ for the \emph{goal function $2$}  by a linear-fractional map while the corresponding optimal values are related by Zhukovskii function:
$$
R_2=\frac{1+R_1}{1-R_1}\cdot\frac{1-\lambda^2}{1+\lambda^2};
\qquad \frac1\mu=\frac12(\lambda+\frac1\lambda),
\qquad 0<\lambda,\mu<1.
$$

It is also easy to notice that two band sets $E$ related by a 1D real projective transformation (i.e. linear fractional one)
have solutions related by pre-composition with the same transformation. Therefore the difference between the corresponding problems is not essential. For this reasons we consider the filter optimization problems up to 
the left action of $PSL_2(R)$ on the band sets $E$.

In a short paper \cite{Sti61} Edward Stiefel gives a soluton to three band filter approximation problem 
in the form of equality of two indefinite integrals: the elliptic and genus $g=2$ abelian one
which implies a functional dependence of their upper limits. Under certain conditions on the integrals this functional dependence is the solution to the above uniform rational approximation problem. 

For the technical ease we use the second goal function in what follows, the uniform norm of the rational approximation of the indicator function $S_E(x)$. The goal of our paper is to develop an efficient function-theoretic algorithm for calculation of  Stiefel fractions adapted to problem setting \eqref{Analogue}. 
Looking ahead we present a graph of a Stiefel fraction right now which helps to grasp the essence of the filter optimization problem. More examples are given in the final section of the paper.

\begin{figure}[h]
\centerline{
\includegraphics[scale=.15]{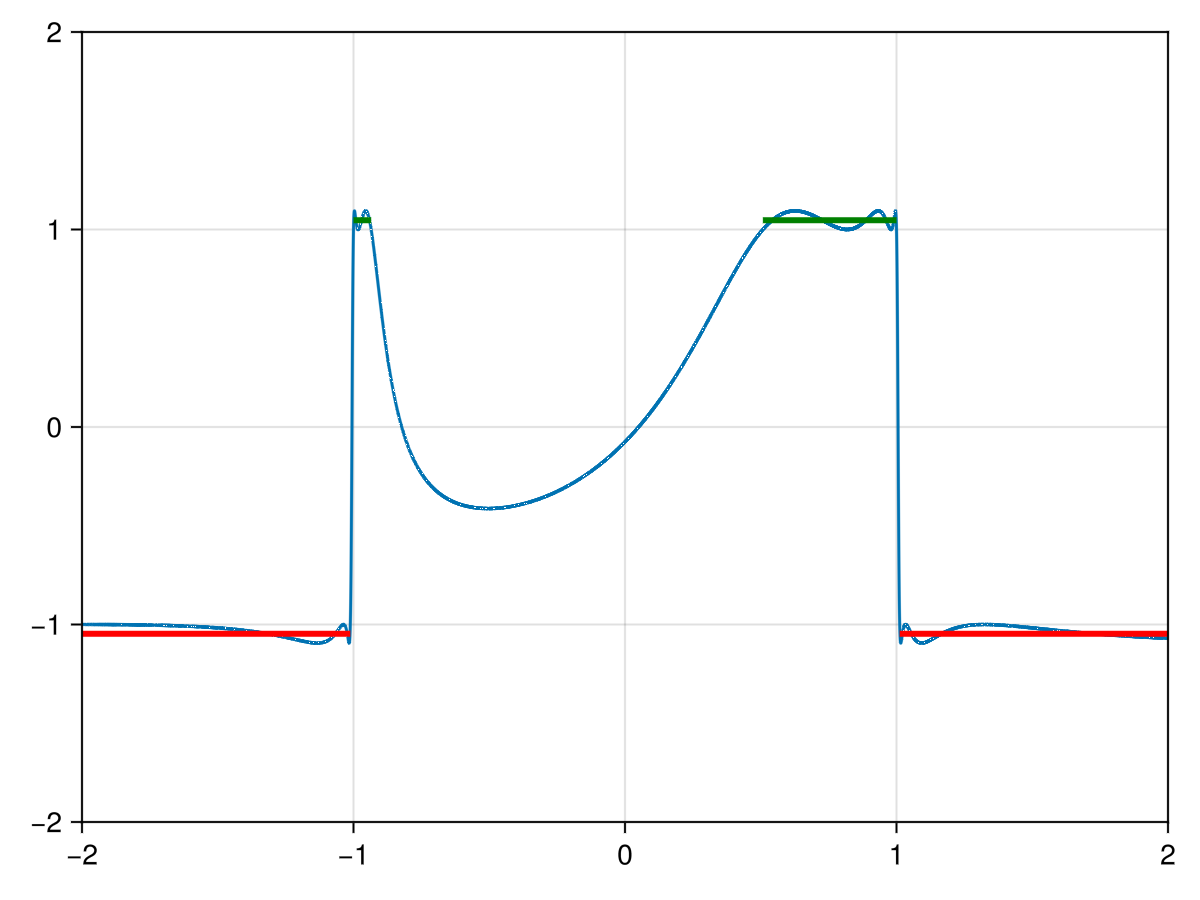}
}
\caption{The graph of a Stiefel fraction of degree $n=10$}
\label{Sti4+4}
\end{figure}

Also we visualise the portion of three band problems solved by 
Stiefel fractions. Recently other solutions to the three band approximation problems were found \cite{B24}.
All of them are given in terms of Chebyshev ansatz \cite{B10,B24} related to Riemann surfaces of genus $g=1$ 
(Zolotarev fraction), $g=2$ (Stiefel fraction) and $g=3$ (some new functions).

{\bf Acknowledgements} This paper was completed when the authors visited Sirius International Mathematical Center
Center. We thank Sirius Center for the hospitality and excellent work conditions. We are also grateful to 
Evgenii Ladeischikov (Steklov Inst.) who has written the first version of computer script for calculation of Stiefel fractions.

\section{Geometric version of Stiefel's solution}
Any local minimum $R(x)$ of filter approximation problem obeys the so called alternation property 
(equiripple in terms of engineers)
\cite {Ach1, Amer, B21}. For the second goal function \eqref{Analogue} its formularion is a bit simpler
than for the first one: the approximation error $\delta(x) := R(x) - S_E(x)$
of  a full degree $n$ minimizer $R(x)$ should have at least $2n + 2$ alternation points $a\in E$, where 
$\delta(a) = \pm||\delta||_C(E)$ with consecutive change of sign. This criterion reduces the solution of 
filter aproximatin problem to construction of rational functions with suitable oscilatory behaviour on the given 
active bands of a filter.

\subsection{Geometric construction for Stiefel fraction}
Below we give a geometric interpretation of Stiefel solution which becomes a Zolotarev fraction 
in the certain limit. First we define two transcendential functions which parametrically 
determine a rational dependency.

Let $\Pi(\tau)$, $\tau\in i(0,\infty) $, be a rectangle of dimensions $2\times|\tau|$:
\be
\label{Rectangle}
\Pi(\tau):=\{u\in \mathbb{C}:\quad -1< Re~u< 1,~~ 0< Im~ u< |\tau| \}.
\ee
Let $x(u|\tau)$ be a confomal map of the rectangle $\Pi(\tau)$ to the upper half plane fixing three points 
$u=0, \pm1$ on the boundary. The latter function has a known representation in terms of elliptic sine or standard 
elliptic theta functions \cite{Ach1}:
\be
\label{exteriour}
x(u|\tau)= sn(K(\tau)u|k)=\frac{\theta_3}{\theta_2}~\frac{\theta_1(u/2)}{\theta_0(u/2)},
\ee
where $sn(\cdot|k)$ is the elliptic sine of modulus $k=k(\tau)=(\theta_2/\theta_3)^2$, 
$K(\tau)=\frac\pi2\theta_3^2$ is the complete elliptic integral; and theta without argument is the 
theta constant i.e. the value of theta for the argument $u=0$ and depending on assumed elliptic modulus $\tau$ only.

Consider also the larger rectangle $\Pi(n\tau)$ tiled by copies of smaller rectangle $\Pi(\tau)$ and  cut along the boundary of one of the tiles:
\be
\label{Pi0}
\Pi(\tau; n,m; h):= \Pi(n\tau)\setminus\{[h,1]+m\tau\},
\qquad n=2,3,\dots;
\quad m=1,2,\dots,n-1;
\quad  -1<h<1.
\ee
One can conformally map the slit rectangle \eqref{Pi0} to the open upper half-plane $\HH$ by a function
$x(u|\tau; n, m; h)$ with the argument $u$ and parameters $\tau,\dots, h$. For convenience we do not 
normalize this map, so the function is defined up to left action of $PSL_2(\mathbb{R})$ 
represented by real linear fractional maps conserving the upper half plane. 
This function does not admit a simple formula like \eqref{exteriour}, however its analogue will be elaborated in the next section. The inverse of this function is a Schwarz-Christoffel integral which in our case is a holomorphic abelian integral on a  suitable genus two Riemann surface. This coincidence will be exploited in what follows.

We deliberately use the same name $x(\dots)$ for two different functions. 
It should not lead to a confusion as the functions differ by the number and type of their arguments.

\begin{thrm}\cite{B24}
Two above introduced transcendential functions parametrically define a rational function:
\be
\label{Sti}
St(z|\tau; n,m; h):=x(u|\tau);
\qquad z:=x(u|\tau; n, m; h),
\qquad  u\in \Pi(\tau; n, m; h).
\ee
Here $St(z|\dots)\in\RR(z)$ and $\deg ~St=n$.
\end{thrm}
{\bf Proof sketch.}
Consider the composite map $St(z|\dots):$ $\HH\to\Pi(\tau; n,m; h)\to\hat{\CC} $ arising after excluding the variable $u$ from formulae \eqref{Sti}. 
It is a meromorphic function in the upper half plane which takes purely real values at the boundary (due 
to quantization of the height of the slit). Therefore it has the compulsory analytic continuation to the whole
Riemann sphere. This function has exactly $n$ simple poles corresponding to $u\in\tau(2\ZZ+1)$ 
(two of those glue in the case $h=0$ and odd $m$). Hence, $St(z|\dots)$ is a real rational function of degree $n$. ~~\bl

\begin{rmk}
This phenomenon is close to the parametric representation of classical Chebyshev polynomials:
$T_n(\cos(u))=\cos(nu)$ with more involved transcendential functions instead of trigonometric ones.
\end{rmk}

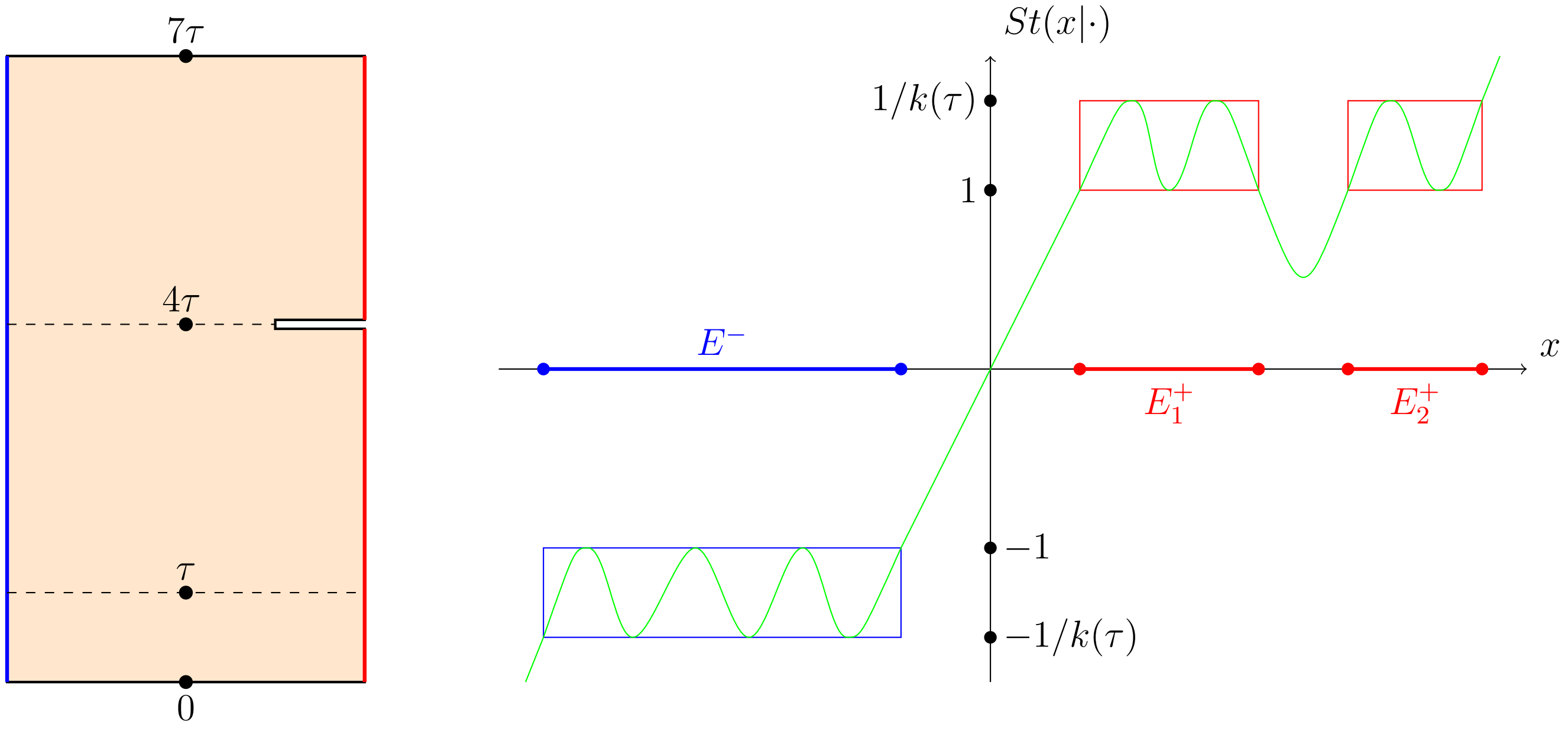
\begin{figure}

\begin{tikzpicture}
\def\L{-11}
\def\R{-7}
\def\midpt{-9}

\draw[thick, fill=orange!20] (\L, -3.5) --
node[midway, below] {$0$}(\R,-3.5) -- 
(\R,0.45) -- (\R-1,0.45) -- (\R-1,0.55)--
(\R,0.55) -- (\R,3.5) --
node[midway, above] 
{$7\tau$}  (\L,3.5) -- cycle;

\draw[very thick, blue] (\L, -3.5) -- (\L, 3.5);
\draw[very thick, red]  (\R,-3.5) --  (\R,0.45);
\draw[very thick, red]  (\R,3.5) --  (\R,0.55);

\draw[dashed, thin] (\L, -2.5) --  (\R,-2.5) 
node[midway, above] {$\tau$};
\draw[dashed, thin] (\L, 0.5) -- (\R-1,0.5)
node[pos=.65, above] 
{$4\tau$};

\foreach \j in {0,1,4,7}
\draw[fill] (\midpt, -3.5+\j) circle (2pt);

\draw[->]    (-5.5,0) -- (6,0)
 node[pos=1.,above right] {$x$};
 \draw[->]    (0,-3.5) -- (0,3.5)
 node[pos=1.,above right] {$St(x|\cdot)$};
 
 \def\stopDistr{{-5,-4.5, -4,-3.3, -2.7,-2.1,-1.6,-1}}; 
 
 \def\passDL{{1, 1.6,2,  2.5, 3}};
 \def\passDR{{4, 4.5, 5, 5.5}};
 \def\devmin{2};
 \def\devmax{3};
 
\draw[blue]  (\stopDistr[0],-\devmax) rectangle 
(\stopDistr[7],-\devmin);
\draw[red]  (\passDL[0],\devmin) rectangle (\passDL[4],\devmax);
\draw[red]  (\passDR[0],\devmin) rectangle (\passDR[3],\devmax);
 
\draw[green] 
(\stopDistr[0],-\devmax) .. controls (\stopDistr[1]-0.15,-\devmin) ..
(\stopDistr[1],-\devmin);

\foreach \j in {1,3,5}
\draw[green] 
(\stopDistr[\j],-\devmin) .. controls (\stopDistr[\j]+0.2,-\devmin) and (\stopDistr[\j+1]-0.2,-\devmax) .. (\stopDistr[\j+1],-\devmax); 

\foreach \j in {2,4}
\draw[green] 
(\stopDistr[\j],-\devmax) .. controls (\stopDistr[\j]+0.2,-\devmax) and  (\stopDistr[\j+1]-0.2,-\devmin) .. (\stopDistr[\j+1,-\devmin); 

\draw[green] 
(\stopDistr[6],-\devmax) .. controls (\stopDistr[6]+0.15,-\devmax) ..
(\stopDistr[7],-\devmin);

\draw[green] 
(\passDL[0],\devmin) .. controls (\passDL[1]-0.15,\devmax) ..
(\passDL[1],\devmax);

\foreach \j in {1}
\draw[green] 
(\passDL[\j],\devmax) .. controls (\passDL[\j]+0.2,\devmax) and (\passDL[\j+1]-0.2,\devmin) .. (\passDL[\j+1],\devmin); 

\foreach \j in {2}
\draw[green] 
(\passDL[\j],\devmin) .. controls (\passDL[\j]+0.2,\devmin) and  (\passDL[\j+1]-0.2,\devmax) .. (\passDL[\j+1,\devmax); 

\draw[green] 
(\passDL[3],\devmax) .. controls (\passDL[3]+0.15,\devmax) ..
(\passDL[4],\devmin);

\draw[green] 
(\passDR[0],\devmin) .. controls (\passDR[1]-0.15,\devmax) ..
(\passDR[1],\devmax);

\foreach \j in {1}
\draw[green] 
(\passDR[\j],\devmax) .. controls (\passDR[\j]+0.2,\devmax) and (\passDR[\j+1]-0.2,\devmin) .. (\passDR[\j+1],\devmin); 

\draw[green] 
(\passDR[2],\devmin) .. controls (\passDR[2]+0.15,\devmin) ..
(\passDR[3],\devmax);

\def\midpt{3.5}

\draw[green] 
(\passDL[4],\devmin) .. controls 
(\midpt,\devmin-1.3) .. (\passDR[0],\devmin);
\draw[green] (\stopDistr[7],-\devmin) -- (\passDL[0],\devmin); 
\draw[green] (\stopDistr[0],-\devmax) -- +(-0.2,-0.5); 
\draw[green] (\passDR[3],\devmax) -- +(0.2,0.5); 

\fill[] (0,-\devmax) circle (2pt)
node[right] {$-1/k(\tau)$};
\fill[] (0,\devmax) circle (2pt)
node[left] {$1/k(\tau)$};
\fill[] (0,-\devmin) circle (2pt)
node[right] {$-1$};
\fill[] (0,\devmin) circle (2pt)
node[left] {$1$};

\foreach \y in {0,4}
\fill[red] (\passDL[\y],0) circle (2pt);
\foreach \y in {0,3}
\fill[red] (\passDR[\y],0) circle (2pt);
\foreach \y in {0,7}
\fill[blue] (\stopDistr[\y],0) circle (2pt);

\draw[blue, very thick] (\stopDistr[0],0) --
node[midway,above] {$E^-$}
(\stopDistr[7],0);
\draw[red, very thick] (\passDL[0],0) -- 
node[midway, below] {$E_1^+$}
(\passDL[4],0);
\draw[red, very thick] (\passDR[0],0) -- 
node[midway, below] {$E_2^+$}
(\passDR[3],0);
 \end{tikzpicture}
\caption{Geometric representation and qualitative graph of a Stiefel fraction for $n=7,~m=4$} 
\label{StiGraph}
\end{figure}

\subsection{Qualitative graph, alternation points and solutions for two problems}
We draw the qualitative graph of the obtained function $St$ on the right panel of Fig \ref{StiGraph}. 
One can observe that after multiplication by the suitable factor  $2k/(k+1)$, $k=k(\tau)$, 
the function $St(x|\dots)$ acquires many alternation points corresponding to the corners of smaller tiles 
on the picture of the  left panel of Fig. \ref{StiGraph}. Required by alternation theorem \cite{Ach1} 
$2n+2$ alternation points belong to the three band set $E$ defined as follows. 
The stopband $E^-$ is the image of the left side of the slit rectangle \eqref{Pi0} under the map 
$x(u|\tau;n,m;h)$,  $u\in -1+\tau[0,n]$.  The choice of two passbands $E^+$ is twofold: 
they are the image of the right side of the slit rectangle with removed segment of the length not greater that $|\tau|$ and either beginning OR ending at the end point $1+m\tau$ of the slit:
 
\be
 E^+(v):=x(u|\tau;n,m;h), 
 \qquad 
 \begin{array}{ll}
 u\in 1+\tau([0,m] \cup [v,n]), & m\le v\le m+1;\\
 OR\\
u\in 1+\tau([0,v] \cup [m,n]), & m-1\le v\le m.
\end{array}
\label{StiE+}
\ee
Note that one of the passbands depend parametrically of the value $v$ in the specified above limits.

The same is true for the reciprocal function $1/St(x|\dots)$ which also displays the equiripple behaviour on the same band set $E$. We have thus proved 
\begin{thrm}\cite{B24} Both functions $\displaystyle{\frac{2k}{k+1} St(x|\dots)}$ and $\displaystyle{\frac2{k+1} \frac1{St(x|\dots)}}$ are local solutions of the filter approximation problem with given above band set $E$ and the same deviation $\mu=\frac{1-k}{1+k}<1$, $k=k(\tau)$.
\end{thrm} 

\subsection{Topological classes of solutions}
It was known already to the pioneers in this field \cite {Zol77, Ach2} that filter optimization problems
may have several local minima. To avoid the nonuniqueness of the solution E.Stiefel and his pupil R.A.Amer 
\cite{Amer} have suggested to restrict the area of optimization to the functions with certain topological 
properties. Those topological classes in the space of given degree rational functions have several definitions
\cite{Amer, B21}, here we use the following: the parity of the number of zeros of $R(x)$ at each gap 
in the bandset $E$. This is an invariant of the continuous 
deformation of  rational function $R(x)$ provided $\mu[R]<1$ and the degree $n:=\deg R$ is fixed.
Topological class is the part of the input data of the problem along with the bandset $E$ and the degree $n$
which guarantees the uniquness of minimum in the closure of the class \cite{Amer, B21}.

For the three band problem the entries of the topological class $\sigma=(\sigma_1,\sigma_{12},\sigma_2)\in \ZZ_2^3$ of the competing function  correspond  to the gaps between the bands $E^-$ and $E_1^+$; $E_1^+$ and $E_2^+$; $E_2^+$ and $E^-$ respectively. Naturally, only two entries of the tuple $\sigma$ are independent:  $\sigma_1+\sigma_{12}+\sigma_2=\deg(R) ~mod ~2$.
An easy exercise is to check that for the (renormalized) Stiefel function $St(x|\dots)$ the topological class is $(1,0,n+1)~ mod~2$ and for the reverse one $1/St$ it is $(0,0,n)~ mod~2$. Two problems with classes $(*,1,*)$ cannot be solved in terms of the Stifel functions at all. The corresponding solutions were found in \cite{B24} and are related to genus $g=3$ curves.

\subsection{Problem: visualization of solution' modes}
We have just noticed that not every three band problem may be solved in terms of Stiefel fractions.
It is natural to ask what portion of all three band problems one can solve in terms of those functions?
We can start answering the question  comparing the dimension of the set of problems and the degrees of 
freedom in Stiefel's fraction. First, the set of problems is three dimensional: any band set $E$ 
is defined by its six end points, but the sets $E$ related by a real projective transformation 
have essentially the same solution. It is exactly for 
this reason that we do not normalize the conformal mapping of the slit rectangle to the upper half plane.
The number of continuous degrees of freedom in Stiefel's construction is also three: two of them ($\tau$, $h$) 
are inherited from the parameters of $St(\cdot|\tau;\dots;h)$. Yet another is hidden 
in the freedom of choice for the passbands described by the variable $v$ in \eqref{StiE+}. 
So we can assume (this can be proved strictly) 
that locally Stiefel fractions indeed solve all three band filter approximation 
problems. Globally it is not correct \cite{B24}.  It is one of the goals of this paper to visualise 
the domain of the problems space with the solution being a Stiefel fraction.

\section{Calculaton of Stiefel fraction}
Stiefel fraction \eqref{Sti} is a composition of two functions: the exteriour function \eqref{exteriour} 
is the elliptic sine \cite{Ach1}, its computation is well elaborated and the appropriate  computer script is easily available. The interiour function is the conformal mapping of the 
upper half plane to the slit rectangle \eqref{Pi0}. It admits the Christoffel-Schwarz 
integral representation, which in our setting becomes a genus 2 hyperelliptic integral and its accurate calculation is a bit more tricky. For computation of this map we adapt the technology elaborated in \cite{B12} and based on the usage of Riemann theta functions of two variables. Of course there are alternative approaches for the numerical conformal mapping, see e.g. \cite{DT}. The one we choose guarantees the machine precision for the solution. 

Suppose the right angles of the slit rectangle correspond to six cyclically ordered points $x_1$, $x_2$,\dots,$x_6$,
of the extended real line $\hat{\RR}$ so that the tip of the slit corresponds to a point 
$c$ between $x_4$ and $x_5$ see Fig. \ref{BasisHomologies}. The conformal mapping $\HH\ni x \to u\in\Pi(\tau;n,m;h)$
may be represented as a Schwarz-Christoffel integral:
\be
\label{SC}
u=1+const\int_{x_3}^x\frac{(x-c)dx}{\sqrt{\prod_{j=1}^6(x-x_j)}},
\qquad x=x(u|\tau;n,m;h).
\ee
We notice that the integral is a holomorphic abelian one related to a genus $g=2$ curve.
Effective realization of this map is decomposed into several natural steps: first we lift the upper half plane $\HH$ 
to an algebraic curve $X$, then we embed the curve into its Jacobian by the Abel-Jacobi mapping. The image of this mapping may be efficiently localized by solution of just one equation. Then we notice that the first coordinate function $u_1$ of the Jacobian maps the image of $\HH$ to a rectangle with a slit.
A brief survey of theory of Riemann surfaces/ algebraic curves necessary for implementing the algorithm
is given below. Further reading on the topic includes \cite{FK, GH, Mu}.

\subsection{Algebraic curve and moduli space} 
Given six real points $-\infty<x_1<x_2<\dots<x_5<x_6<\infty$ of the real line,
consider a compact genus two Riemann surface $X$ with the equation
(of its affine part):
\be
\label{X}
y^2=\prod\limits_{s=1}^6(x-x_s),
\qquad  (x,y) \in{\CC}^2.
\ee
 This curve admits a \emph{conformal involution} $J(x,y)=(x,-y)$ with six stationary points $p_s:=(x_s,0)$ and an \emph{anticonformal involution} (\emph{reflection})
$\bar{J}(x,y)=(\bar{x},\bar{y})$. The stationary points set of the latter map  has three 
components known as \emph{real ovals} of the curve. Each real oval is an embedded circle \cite{Nat} and doubly covers exactly one of the segments  $[x_2,x_3]$, $[x_4,x_5]$, $[x_6,x_1]\ni\infty$ of the extended real line $\hat{\mathbb R}$. 
We denote those ovals as \emph{first, second and third} respectively.
The  lift of the complimentary set of intervals
$[x_1,x_2]$, $[x_3,x_4]$, $[x_5,x_6]$ to the surface (\ref{X})
 gives us three \emph{coreal ovals} which make up the set of points fixed by another anticonformal involution ${\bar J}J=J{\bar J}$.

\begin{dfn} Moduli space ${\cal M}_2\mathbb{R}_3$ is the space of
genus two Riemann surfaces with reflection and three enumerated real ovals. 
Two surfaces are equivalent iff  there is a conformal $1-1$ mapping between them commuting 
with the reflections and respecting the marking of real ovals.
\end{dfn}

Saying differently \cite{B12}, the moduli space  ${\cal M}_2\mathbb{R}_3$ is the space  of   the cyclically ordered sextuples  of  points $(x_1, \dots, x_6)$ from the extended real line
$\hat{\mathbb{R}}$ modulo the action of $PSL_2(\mathbb{R})$ 
(= real projective transformations conserving the orientation). Topologically it is a 3-cell.

\subsection{Homologies, differentials, periods}\label{SectBasisCycles}
\begin{figure}
\centerline{\includegraphics[scale=1.]{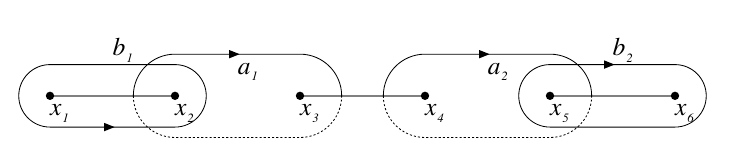}}
\caption{\small Canonical basis in homologies of the curve $X$}
\label{BasisHomologies}
\end{figure}

We fix a special basis in the 1-homology space of the curve $X$ intrinsically related to the latter. 
Two 1-cycles  $a_1$ and $a_2$  are the first and second real ovals shown in Fig. \ref{BasisHomologies}; 
two complimentary cycles $b_1$ and $b_2$ are the coreal ovals of the curve  shown on the same figure and 
oriented so that the intersection matrix takes the canonical form.

The reflection $\bar{J}$ of the surface conserves the $a-$ cycles and reverses the $b-$ cycles.
The basis of holomorphic differentials $du_1$, $du_2$ dual to the basis of cycles
\be
\int_{a_s}du_j:=\delta_{sj};
\qquad s,j =1,2,
\ee
determines Riemann period matrix $\Pi$ with the elements
\be
\label{periods}
\Pi_{sj}:=\int_{b_s}du_j;
\qquad s,j =1,2.
\ee
It is a classical result due to B.Riemann that $\Pi$ is symmetric and has positive definite imaginary part \cite{FK}.
In our case both basic differentials $du_1$, $du_2$ are real (in the usual sense $\bar{J}^*du=\overline{du}$)
because of  the mentioned above symmetry of the homology basis with respect to the reflection. 
In particular this means that the period matrix is purely imaginary: $\Pi=i\Omega$ with real symmetric and positiley definite matrix $\Omega$.

\begin{thrm} E.g. \cite{B12}
The period mapping $\Omega(X)$  is real analytic diffeomorphism from the moduli space ${\cal M}_2\mathbb{R}_3$  to the trihedral cone
\be
\label{PeriodsCone}
0<\Omega_{12}<min(\Omega_{11}, \Omega_{22}).
\ee
\end{thrm}

\subsection{ Jacobian and Abel-Jacobi map }
\begin{dfn}
Given a Riemann period matrix $\Pi$, we define the full rank lattice
\be
\label{Lattice}
L(\Pi)=\mathbb{Z}^2+ \Pi\mathbb{Z}^2 = \int_{H_1(X,\mathbb{Z})}d{\bm u},
\qquad d{\bm u}:=(du_1,du_2)^t,
\ee
in $\mathbb{C}^2$ and the 4-torus $Jac(X):=\mathbb{C}^2/L(\Pi)$ known as a Jacobian of the curve $X$.
\end{dfn}
This definition depends on the choice of the symplectic basis in the lattice $H_1(X,\mathbb{Z})$, 
other choices bring us to isomorphic tori.

It is convenient to represent a point  ${\bm u}\in\mathbb{C}^2$ as a theta characteristic $[\epsilon,\epsilon']$, i.e. a couple of real 2-vectors (columns) $\epsilon, \epsilon'$:
\be
{\bm u}=\frac12(\Pi\epsilon+\epsilon').
\ee
The points of Jacobian  $Jac(X)$ in this notation correspond to
two vectors with real entries modulo 2. Second order points of Jacobian are
represented as $2\times 2$ matrices with ${\mathbb Z}_2$ entries.
Note that we use a nonstandard notation of theta characteristic as two column vectors 
written one after another. Sometimes the transposed matrix is used.

\begin{dfn}
Abel-Jacobi (briefly: AJ) map is a correctly defined mapping from the surface $X$ to its Jacobian.
\be
\label{AJmap}
{\bm u}(p):=\int_{p_1}^p d{\bm u}~~ mod~L(\Pi),
\qquad p_1:=(x_1,0); \quad d{\bm u}:=(du_1,du_2)^t.
\ee
\end{dfn}

From Riemann-Roch formula it follows \cite{FK} that Abel-Jacobi map is a holomorphic embedding 
of the curve into its Jacobian. In Sect. \ref{theta} we give an explicit formula for the image of this embedding for the genus $g=2$ case. Meanwhile we compute the images of the branching points $p_s=(x_s,0)$, $s=1,\dots,6$ of the curve $X$:

\centerline{
\begin{tabular}{c|c|c}
$p$ & ${\bm u}(p)~ mod~ L(\Pi)$ & $[\epsilon, \epsilon']({\bm u}(p))$\\
\hline
$p_2$ & $\Pi^1/2$ & $\tiny
\left[\begin{array}{c} 10\\00\end{array}\right]$\\
$p_3$&$(\Pi^1+E^1)/2$&$\tiny
\left[\begin{array}{c} 11\\00\end{array}\right]$\\
$p_4$&$(\Pi^2+E^1)/2$&$\tiny
\left[\begin{array}{c} 01\\10\end{array}\right]$\\
$p_5$&$(\Pi^2+E^1+E^2)/2$&$\tiny
\left[\begin{array}{c}01\\11\end{array}\right]$\\
$p_6$&$(E^1+E^2)/2$&$\tiny
\left[\begin{array}{c}01\\01\end{array}\right]$\\
\label{AJPj}
\end{tabular}}
\noindent
here $\Pi^s$ and $E^s$ are the $s$-th columns of the period and the identity matrix respectively. 

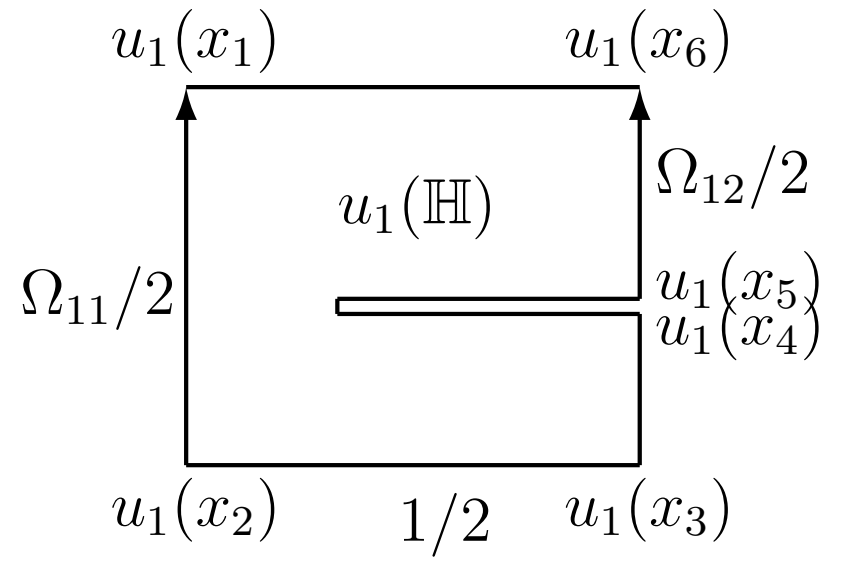
\begin{figure}[h!]
\begin{picture}(100,40)(-60,0)
\thicklines
\put(10,5){\vector(0,1){25}}
\put(10,5){\line(1,0){30}}
\put(40,5){\line(0,1){10}}
\put(40,15){\line(-1,0){20}}
\put(20,15){\line(0,1){1}}
\put(20,16){\line(1,0){20}}
\put(40,16){\vector(0,1){14}}
\put(40,30){\line(-1,0){30}}
\put(5,1){$u_1(x_2)$}
\put(35,1){$u_1(x_3)$}
\put(41,13){$u_1(x_4)$}
\put(41,16){$u_1(x_5)$}
\put(35,32){$u_1(x_6)$}
\put(5,32){$u_1(x_1)$}
\put(-1,15){$\Omega_{11}/2$}
\put(41,23){$\Omega_{12}/2$}
\put(24,0){$1/2$}
\put(20,21){$u_1(\mathbb{H})$}

\end{picture}
\caption{\small The image of the upper half plane $u_1(\mathbb{H})$ under the first coordinate 
function of the Jacobian}
\label{URect}
\end{figure}

The basic differential $du_1$ on $X$ is positive \cite{Nat} on the first real oval, 
has two zeros interchanged by the involution $J$ on the second real oval. It is not difficult to check (see details in \cite{B12}) that the first jacobi coordinate function $u_1(p)$ maps a suitable lift of the 
upper half plane $\HH\subset X$ to  a rectangle with a slit shown in Fig. \ref{URect}.
The boundary orientation of the lifted half plane should coinside with the orientation of the first real oval at the segment of their intersection. This is exactly what we need for the calculation of Stiefel fraction.  In out case 
\be
\label{interiour}
u=1+4\int_{x_3}^x du_1,
\qquad x=x(u|\tau,n,m;h)\in\HH.
\ee
The dimensions of the polygon $u_1(\HH)$ are related to the entries of the
period matrix. Therefore in our case we know two of three independent entries:
\be
\label{Omega}
2i\Omega_{11}=n\tau;
\qquad
2i\Omega_{12}=(n-m)\tau;
\ee
and the length of the slit is a function of the remaining unknown  element $\Omega_{22}$ of the 
period matrix.

\subsection{Theta functions on genus two surface}\label{theta}
Two auxuliary problems related to conformal mapping of the slit rectangle
will be effectively solved in terms of Riemann theta functions:
\begin{itemize}
\item Localization of the curve inside its Jacobian;
\item Representation of projection $p=(x,y)\to x$ of the curve $X$ to the sphere;
\end{itemize}

\begin{dfn}
Let ${\bm u}\in\mathbb{C}^2$ and $\Pi\in\mathbb{C}^{2\times2}$ be a Riemann
matrix, i.e. $\Pi=\Pi^t$ and $Im~\Pi>0$. The theta function of those two arguments is the following Fourier series

\be
\label{thetadef}
\theta({\bm u}, \Pi):=\sum\limits_{m\in\mathbb{Z}^2}
\exp(2\pi im^t{\bm u}+\pi i m^t\Pi m).
\ee
They also consider theta function with characteristics
which is the slight modification of the above theta.
\be
\label{thetachardef}
\theta[2\epsilon, 2\epsilon']({\bm u}, \Pi):=\sum\limits_{m\in\mathbb{Z}^2}
\exp(2\pi i(m+\epsilon)^t({\bm u}+\epsilon')+\pi i (m+\epsilon)^t\Pi (m+\epsilon))
\ee
$$
=
\exp(i\pi\epsilon^t\Pi\epsilon+2i\pi\epsilon^t({\bm u}+\epsilon'))
\theta({\bm u}+\Pi\epsilon+\epsilon',\Pi),
\qquad \epsilon,\epsilon'\in\mathbb{R}^2.
$$
The matrix argument $\Pi$ of theta function is usually omitted when it is clear which matrix is meant.  
Omitted vector argument $u$ is supposed to be zero and the appropriate function of $\Pi$ is called 
the theta constant:
$$
\theta[\epsilon, \epsilon']:=\theta[\epsilon, \epsilon'](0, \Pi).
$$
\end{dfn}
The convergency of these seria grounds on the positive determinacy of
$Im~\Pi$.  The series has high convergency rate with well controlled accuracy \cite{DHB}.

\begin{rmk}\label{RemTheta}
(i) Theta function with integer characteristics
$[\epsilon, \epsilon']$ is either even or odd depending on the parity of the inner product $\epsilon^t\cdot\epsilon'$. In particular, all odd theta constants vanish.



(ii) It is convenient to represent integer theta characteristics as the sums of AJ images of the branchpoints, keeping only the indexes of those points:
$[sk..l]$ means the sum modulo 2 of the theta characteristics of points
$u(p_s)$, $u(p_k)$, $\dots,u(p_l)$, e.g. $[35]$ in this notation is the so called vector of Riemann constants represented by theta characteristics.
\end{rmk}

\subsection{Image of Abel-Jacobi map}
The location of a genus 2 curve embedded to its Jacobian by AJ map may be reconstructed by solving a single equation.

\begin{thrm}[Riemann]
A point ${\bm u}$ of Jacobian lies in the image ${\bm u}(X)$ of Abel-Jacobi map
if and only if
$$
\theta{\tiny\left[\begin{array}{c}10\\11\end{array}\right]}({\bm u})=0.
$$
\end{thrm}
For short and easy proofs of this and other theorems of this section see e.g. \cite{B12}.

\subsection{Projection to the sphere}\label{projection}
Any meromorphic function on the curve may $X$ be effectively calculated via the
Riemann theta functions once we know its zeros and poles. Take for instance the degree 2 function $x$ on the hyperelliptic curve (\ref{X}). This projection is unique if normalized  e.g. as follows:
$x(p_s)=0$, $x(p_j)=1$, $x(p_l)=\infty$, $s,j,l$ is a positive triple from the index set 
$\{1,2,\dots,6\}$.

\begin{thrm}\cite{FK,B12}
 The hyperelliptic projection $x(p)$ with the above normalization 
 in terms of jacobian coordinates has the appearance:
\be
\label{xofP}
x(p)= \pm\frac{\theta^2[lkj35]}{\theta^2[skj35]}
      \frac{\theta^2[sk35]({\bm u}(p))}{\theta^2[lk35]({\bm u}(p))},
\qquad k\neq s,l,j,
\ee
where the sign $\pm$ in the latter formula
depends on the parity ($+$ even /$-$ odd)
of the scalar product $\epsilon(j)\cdot(\epsilon'(s)+\epsilon'(l))$, where
$[\epsilon(*),\epsilon'(*)]$ corresponds to a branchpoint ${\bm u}(p_*)$ as in Tab. \ref{AJPj}.
\end{thrm}

\begin{xmpl}
For instance the normalized projection with
$x_1=0$, $x_2=1$, $x_6=\infty$, 
is given by the formula
\be
\label{x136}
x(p)= -\frac{\theta^2[1]}{\theta^2[6]}
      \frac{\theta^2[26]({\bm u}(p))}{\theta^2[2]({\bm u}(p))},
\ee
All other normalizations are obtained by the left action of $PSL_2(\RR)$ on this map.
\end{xmpl}

\subsection{Assembly of the algorithm}
The input data of the algorithm are $\tau\in i(0,+\infty)$, integers $n\ge2$, $m=1,2\dots,n-1$ and instead of $h$ we use the functionally dependent value $\Omega_{22}>(n-m)|\tau|$. When 
$\Omega_{22}$  changes from 
$\Omega_{12}$ to $+\infty$, the value $h$ governing the length of the slot monotonically rises 
from $-1$ to $+1$, see e.g. \cite{B12}.

We initialize the period matrix $\Pi=i\Omega$ by \eqref{Omega} and the value of $\Omega_{22}$.
Then, given the argument $x\in\HH$ of the Stiefel fraction, we calculate a solution ${\bm u}=(u_1,u_2)$ of the $2\times 2$ transcendential system 
\be
\theta\left[
{\tiny
\begin{array}{c} 10\\11\end{array}}
\right]
({\bm u})=0;  
\qquad -\frac{\theta^2\left[
{\tiny
\begin{array}{c} 00\\00\end{array}}
\right]}
{\theta^2\left[
{\tiny
\begin{array}{c} 01\\01\end{array}}
\right]}
\frac{\theta^2\left[
{\tiny
\begin{array}{c} 11\\01\end{array}}
\right]
({\bm u}(p))}{\theta^2
\left[
{\tiny
\begin{array}{c} 10\\00\end{array}}
\right]({\bm u}(p))}
      =x,
\qquad k\neq 1,3,6.      
\ee
There are exactly two solutions on this system modulo the periods lattice $L(\Pi)$ which differ by the sign and correspond to two different lifts of the upper half plane to the surface $X$. We choose the solution $\bm u$ with theta characteristics belonging to the block $\displaystyle\left[
\begin{array}{c} -I~I\\-I~I\end{array}\right] ~ mod~2$, $I:=(0,1)$ of the Jacobian, see explanation in \cite{B12}.

Once the auxiliary jacobian coordinates are found, the value of $St(x|\dots)$ is equal to the value of function \eqref{exteriour} with the argument 
$$
u=4u_1+(n-m)\tau-1.
$$
For the chosen normalization of the projection $x(p)$
the stopband $E^-:=[0,1]$, and the right end of $E_2^+$ equals $\infty$. Three other ends of the bandset are ratios of theta constants obtained by substituting the AJ images of the branchpoints to the formula \eqref{x136}.  Yet unknown endpoints of the  maximal pass bands $E_1^+$, $E_2^+$ (when $v=m$) take the appearance 
\be
x_3=\frac{
\theta^2\left[
{\tiny
\begin{array}{c} 00\\00\end{array}}
\right]
\theta^2\left[
{\tiny
\begin{array}{c} 00\\01\end{array}}
\right]
}
{\theta^2\left[
{\tiny
\begin{array}{c} 01\\00\end{array}}
\right]
\theta^2\left[
{\tiny
\begin{array}{c} 01\\01\end{array}}
\right]
},
\quad
x_4=\frac{
\theta^2\left[
{\tiny
\begin{array}{c} 00\\01\end{array}}
\right]
\theta^2\left[
{\tiny
\begin{array}{c} 10\\10\end{array}}
\right]
}
{
\theta^2\left[
{\tiny
\begin{array}{c} 01\\00\end{array}}
\right]
\theta^2\left[
{\tiny
\begin{array}{c} 11\\11\end{array}}
\right]
},
\quad
x_5=\frac{
\theta^2\left[
{\tiny
\begin{array}{c} 00\\00\end{array}}
\right]
\theta^2\left[
{\tiny
\begin{array}{c} 10\\10\end{array}}
\right]
}
{
\theta^2\left[
{\tiny
\begin{array}{c} 11\\11\end{array}}
\right]
\theta^2\left[
{\tiny
\begin{array}{c} 01\\01\end{array}}
\right]
}.
\label{branch1}
\ee
respectively.

\section{Numerical experiments}
Below we give two applications of the elaborated algorithm for calculation of Stiefel functions and visualisation of three band problems admitting Stiefel solutions.

\subsection{Several typical Stiefel fractions}
\begin{figure}[h]
\centerline{
 \includegraphics[scale=.15]{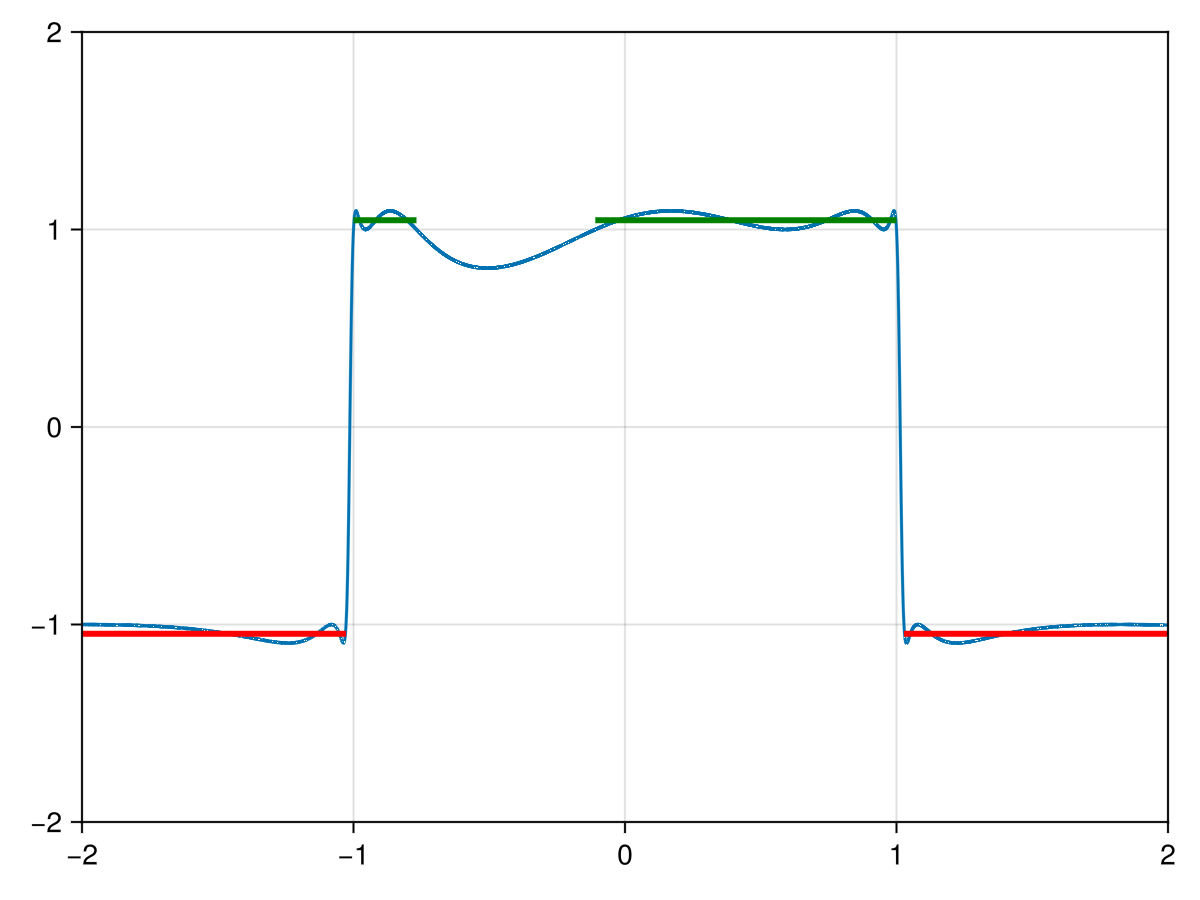}
 \includegraphics[scale=.15]{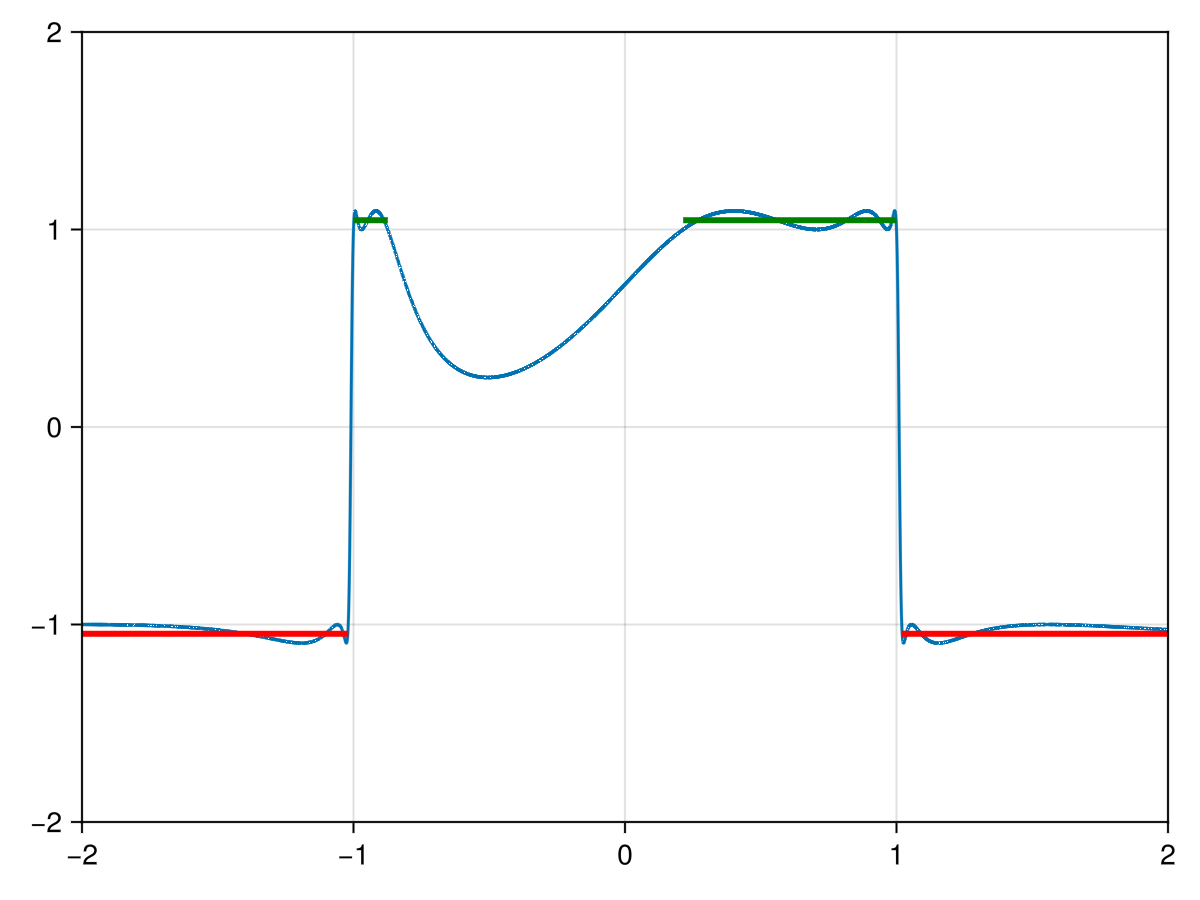}
 \includegraphics[scale=.15]{down_3.png}}
\centerline{
 \includegraphics[scale=.15]{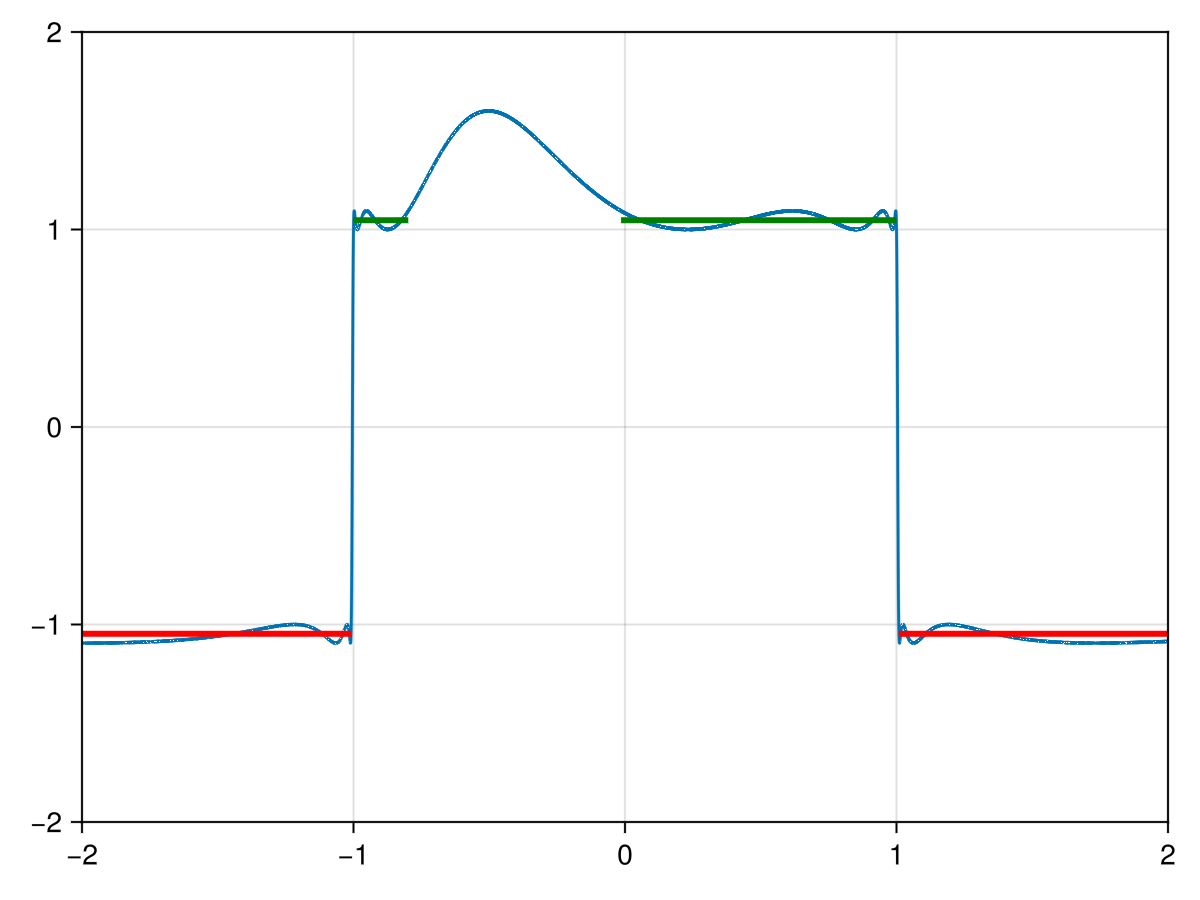}
 \includegraphics[scale=.15]{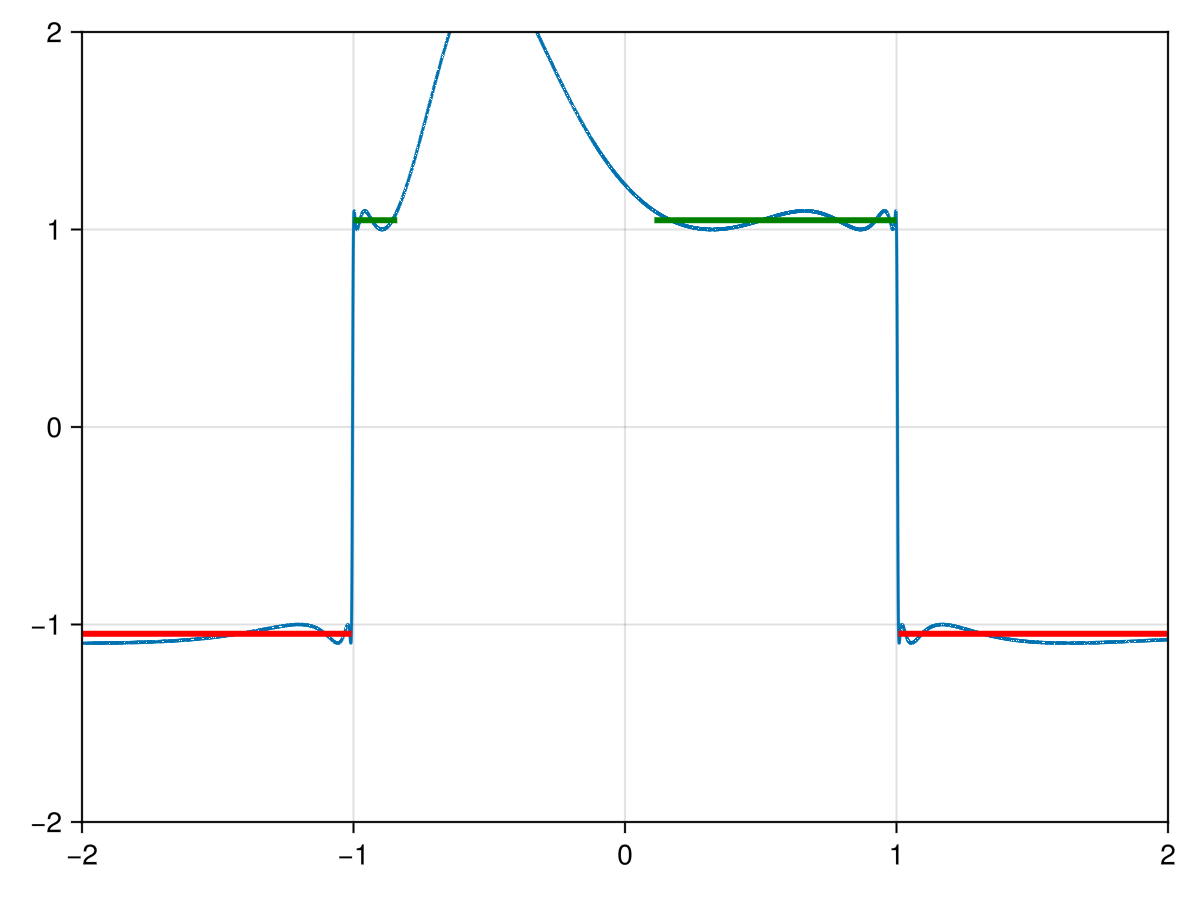}
 \includegraphics[scale=.15]{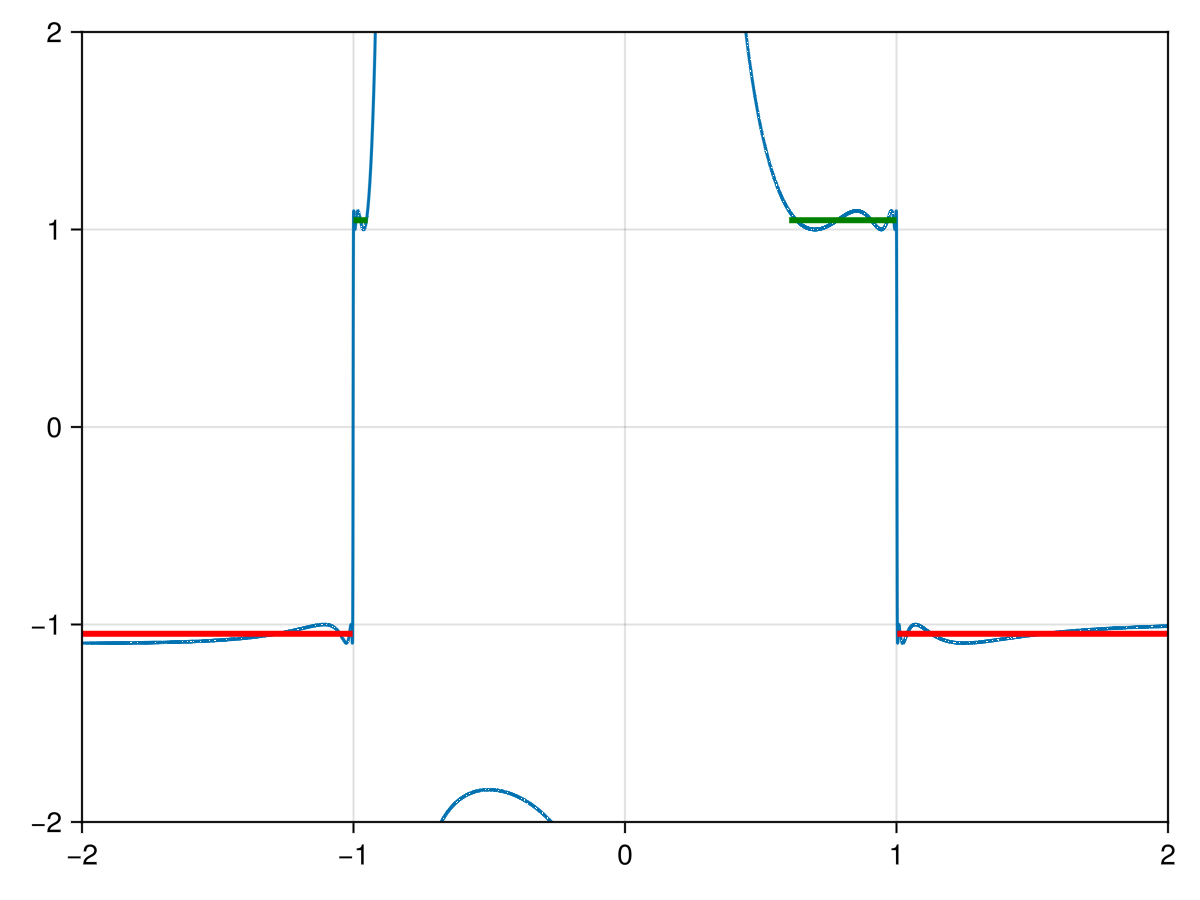}} 
\caption{Stiefel fractions of degree $n=10$ and oscillation number $m=4$ (top row), degree $n=12$ and $m=5$ (bottom row),  $\tau=0.7i$ and different values of $h$}
\label{StiFigs} 
\end{figure}

To compute a Stiefel fraction $St(\dots)$ one can use the positions of its zeros and poles
which correspond to $u\in 2\ZZ\tau$ and $(2\ZZ+1)\tau$ respectively. Then the obtained fraction may be normalized e.g. by  its known value at an endpoint of the bandset $E$. The hump of the graph in the passband gap 
is oriented down for the even oscillation number $m$ or up if it is odd. The hump disappears 
as the parameter $h$ governing the length of the slot tends to $h=1$, the solution itself becomes a Zolotarev fraction and may be parametrically expressed in terms of elliptic functions. In the latter limit case $x(u|\tau;m,n;1)=x(u|n\tau)$ 
(modulo the action of projective transformation).

\subsection{Domains corresponding to Stiefel solutions}
We know that Stiefel fraction $St(x|\dots)$ depends on continuous parameters $\tau$ and $h$ varying in the known limits as well as integer parameters $n>1$ and $0<m<n$.  Also, the fixed function gives a solution of a 1-parametric set of problems,
described in \eqref{StiE+} by parameter $v$ (modulo action of projective group on the sets $E$ as a whole which gives yet three inessential degrees of freedom). The boundary between Stiefel and Zolotarev solutions are given by equality $h=+1$
two other continuous degrees  of freedom $\tau$, $v$ paramertrize this boundary and integer parameters $m,n$ enumerate different components of this boundary. Another boundary between problems with Stiefel solutions and those described by 
genus $g=3$ curves \cite{B24} is given by the equalities $v=m\pm1$ for the fixed $m$. This gives the decomposition  
of the space of problems $[E,\sigma,n]$ into $3n-3$ domains for $\sigma=(1,0,n-1)~mod~2$. More on this decomposition will be given in subsequent publications.

\begin{figure}[h!]
 \centering{
 \includegraphics[scale=.16]{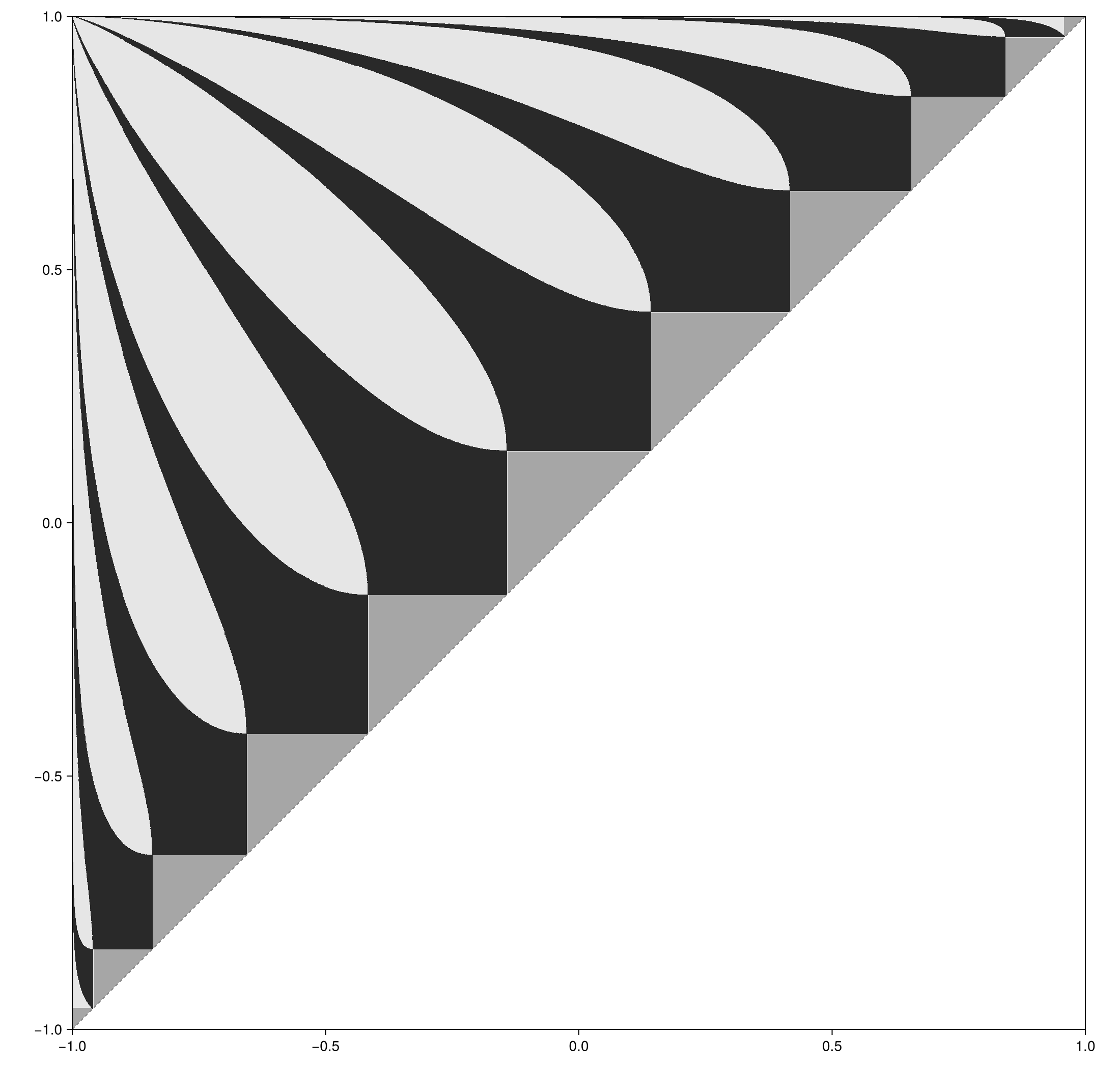}
 \includegraphics[scale=.16]{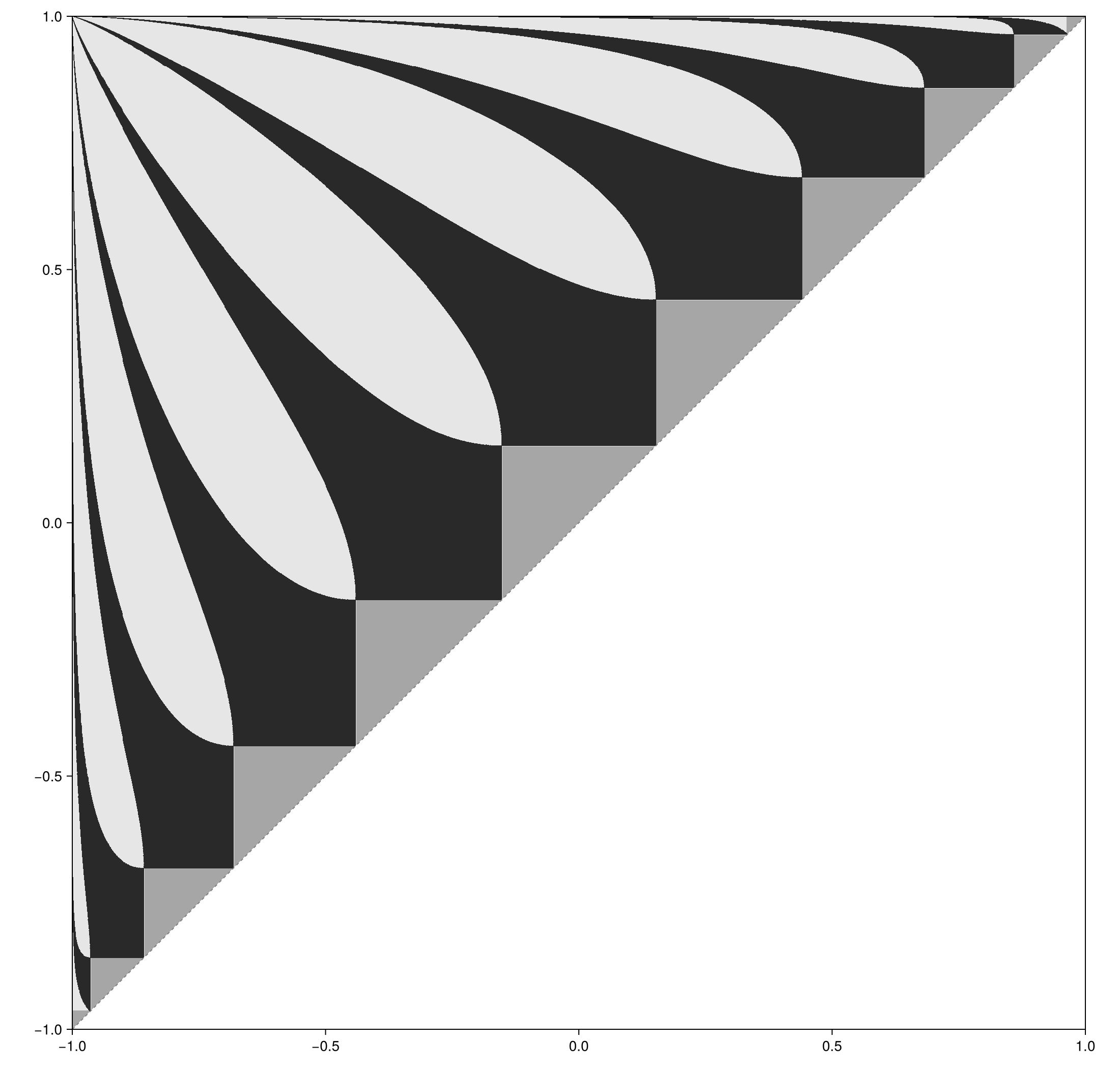}
 \includegraphics[scale=.08]{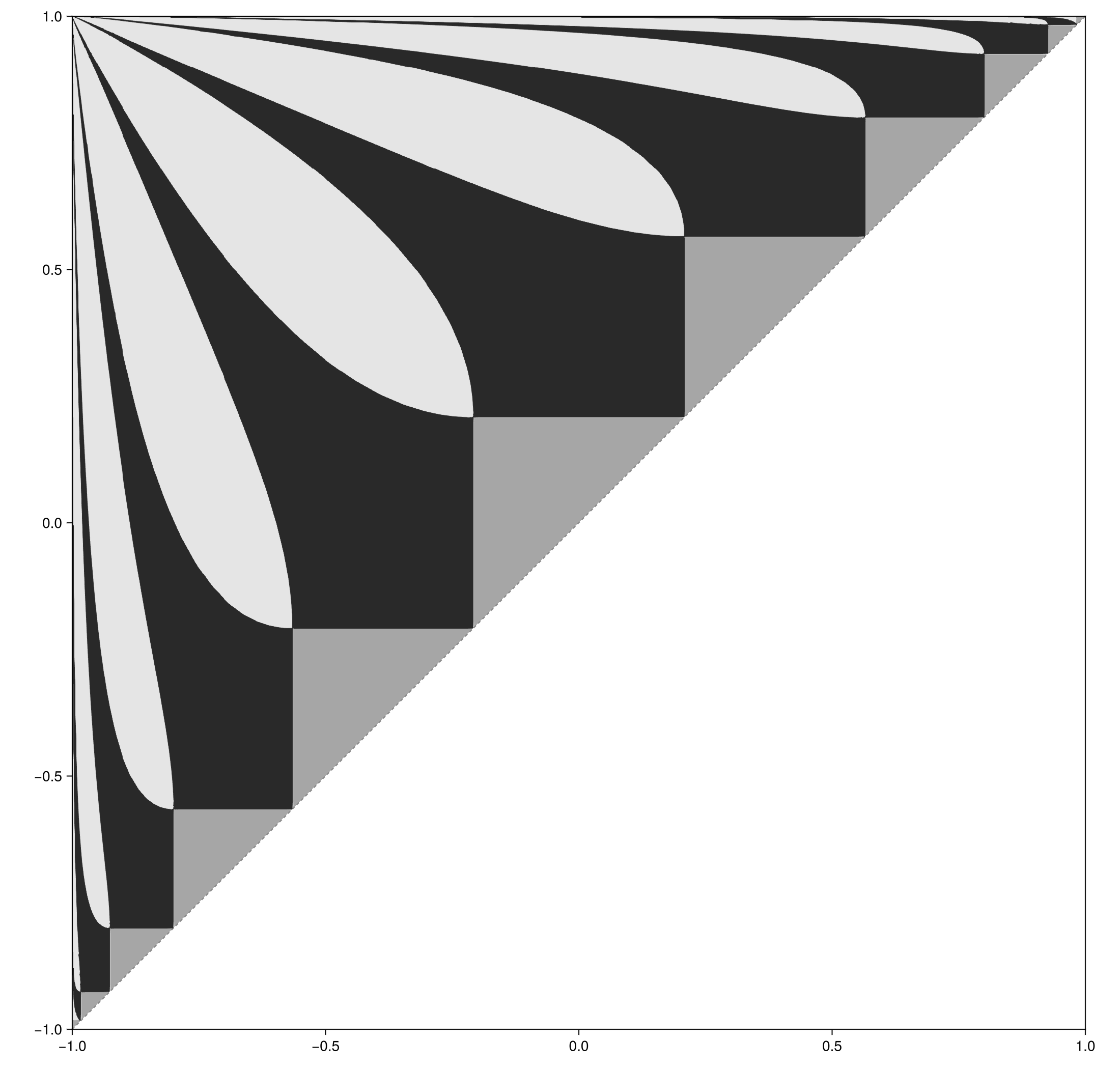}}
 \centering{
 \includegraphics[scale=.08]{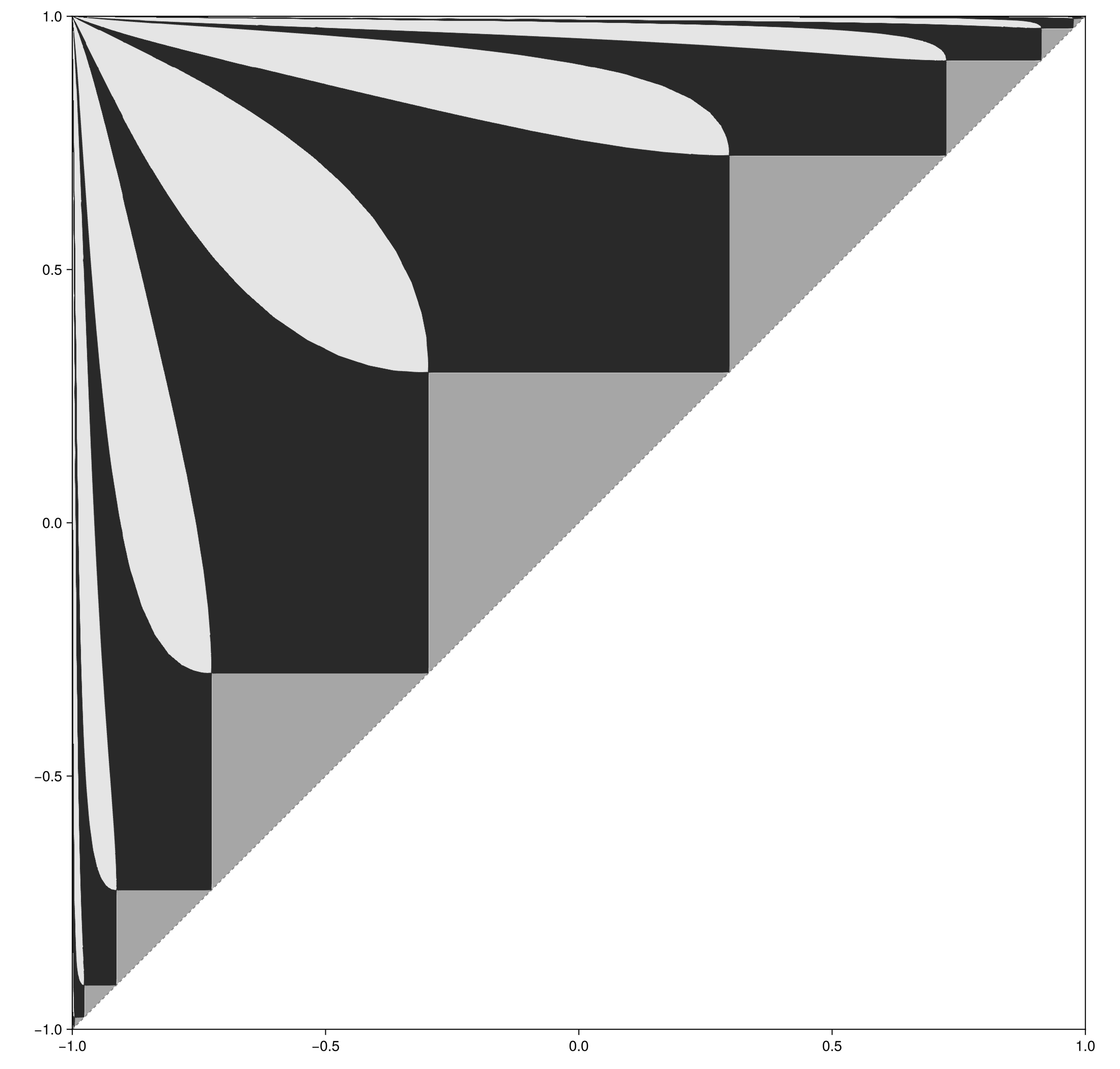}
 \includegraphics[scale=.08]{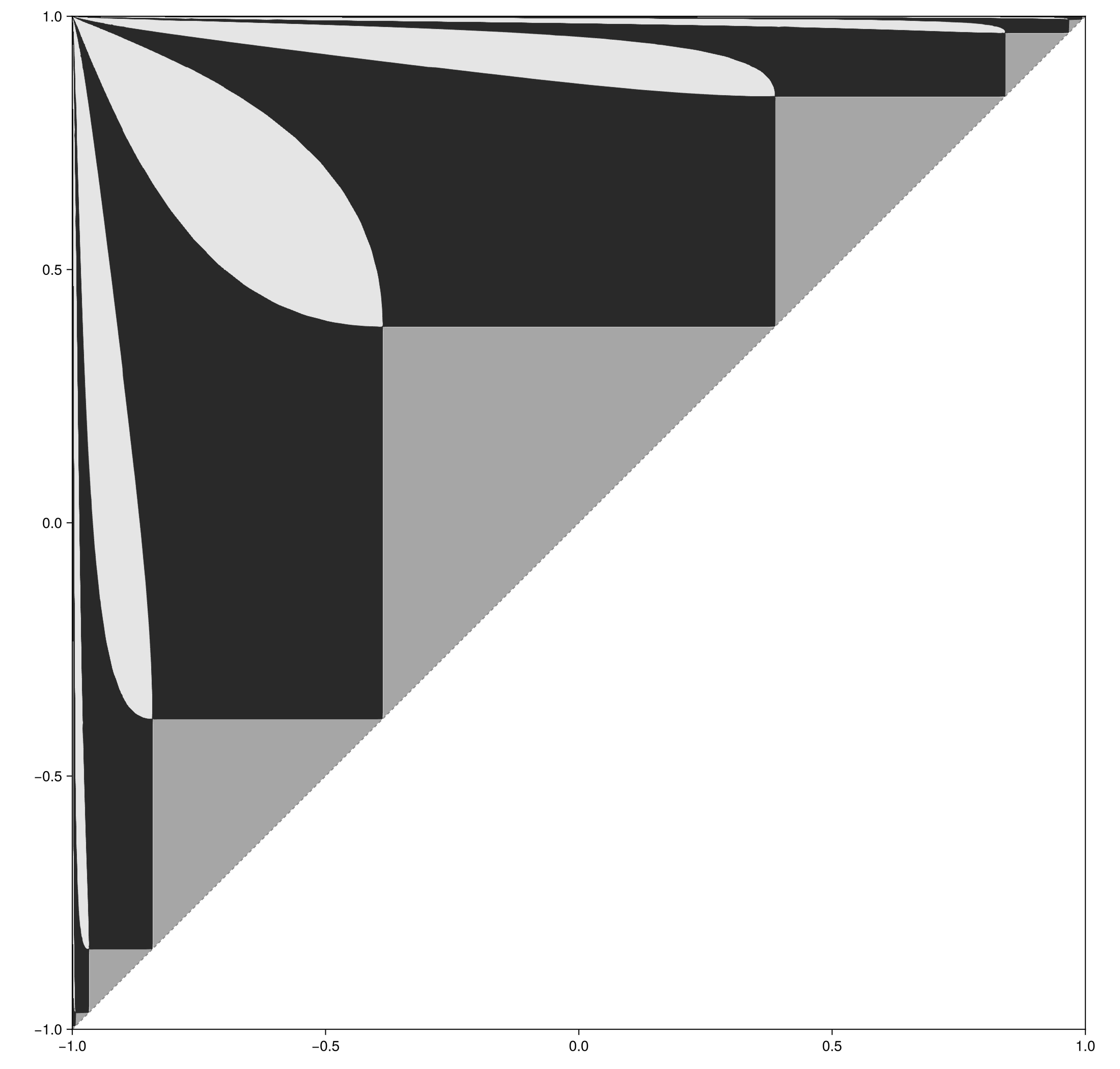}
 \includegraphics[scale=.08]{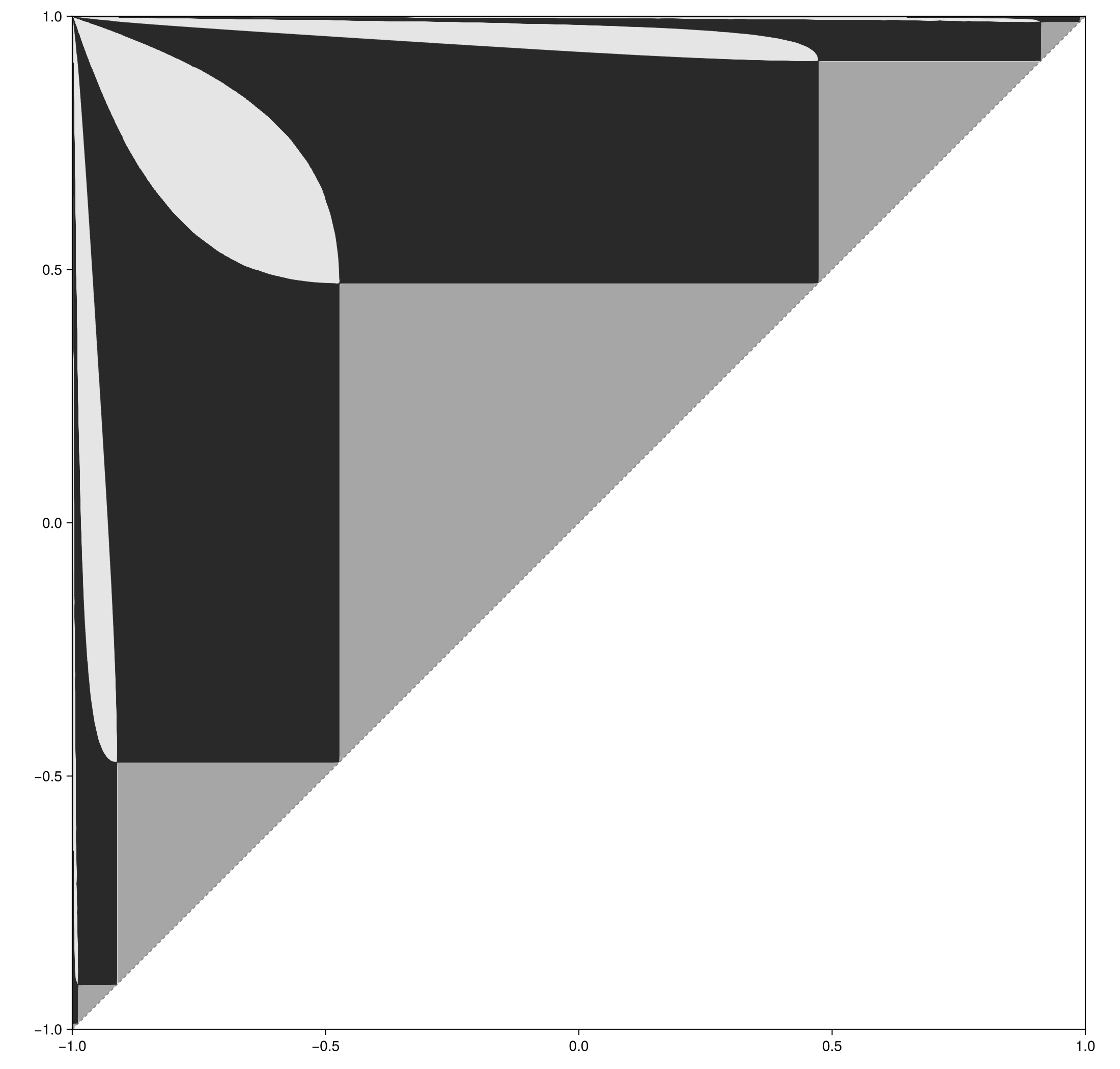}}
\caption{The space of three band problems of degree $n=11$ and topological class $\sigma=(1,0,0)$: several planar sections $1/\gamma= 10,2,1.1$ (top row) and $1/\gamma= 1,1.01,1.001$ (bottom row). The domains corresponding to 
Stiefel, Zoloteraev and genus three solutions are shown in dark,  gray and light gray hue respectively}
\label{StiSaw}
 \end{figure}

For the better visibility of phase diagrams of three band filter problem we use the folowing normalization of the bandset $E$:
\be
\label{NormE}
E^-:=[1/\gamma,-1/\gamma]\ni\infty; 
\qquad E_1^+:=[-1,\alpha]; 
\quad E_2^+:=[\beta,1]. 
\ee
Therefore the band set $E$ acquires the coordinates $(\alpha,\beta,\gamma)$ varying in a triangular prism
$$
\gamma\in (0,1); \qquad -1<\alpha<\beta<1.
$$
We consider the sections of this prism by the planes with constant coordinate $\gamma$. The pictures of the decomposition 
of the space of three band filter problems into Zolotarev, Stiefel and $g=3$ modes are shown on Fig. \ref{StiSaw}.

\begin{verbatim}
G.I.Marchuk Institute for Numerical Math., 
Russian Academy of Sciences, 
Address: Russia 119991 Moscow ul. Gubkina, 8   
Email: a.bogatyrev@inm.ras.ru,   t.garaev@inm.ras.ru  
\end{verbatim}


\begin{thebibliography}{22}
\bibitem{Zol77}
Zolotar\"ev, E.I., Application of Elliptic Functions to Questions on Functions
Deviating Least and Most from Zero, \emph{Zap.\ Imp.\ Akad.\ Nauk St.~Petersburg},
vol.~30, no.~5 (1877), pp.~1--71.
\bibitem{Ach1}
Achieser, N.I., Elements of the theory of elliptic functions, 
(Russian). Leningrad 1948; Transl. of Math. Monographs, 79, AMS, Providence, RI, 1990. viii+237 pp.
\bibitem{Ach2}
Achieser, N.I., Theory of Approximation (Monograph) --M: Nauka, 1965;  Dover Publications, 1992.
\bibitem{Sti61} E. Stiefel,
Le problème d'approximation dans la théorie des filtres électriques,
//Colloque sur l'analyse numérique à Mons, 1961, pp. 81--87. 
\bibitem{Amer}
R. A.-R. Amer and H. R. Schwarz, Contributions to the approximation problem of electrical filters. 
Mitt. Inst. Angew. Math.
(1964), No. 9, Birkhäuser, Basel. 99pp. 
See also R.A.-R. Amer's PHD Thesis, ETH, 1964. 
\bibitem{B10}
Bogatyrev, A.B., Chebysh\"ev Representation of Rational Functions, \emph{Sb.\ Math.}\ (Engl.\ Transl.), 2010, vol.~201, no.~11--12, pp.~1579--1598.
\bibitem{B21} A.B. Bogatyrev, “Projective view at optimization problem for multiband filter”, Proceedings of AMS, 149:7 (2021), 3021-3035
\bibitem{B24} A.B. Bogatyrev, 'Stiefel filters', Trans. Moscow Math Soc., vol 85:2 (2024), MCCME, pp. 323--337   
\bibitem{B12} A.B. Bogatyrev, “Conformal mapping of rectangular heptagons”, Sb. Math., 203:12 (2012), 1715--1735;
\bibitem{BG18} A.Bogatyrev and O.Grigor'ev, “Conformal mapping of rectangular heptagons II”, CMFT, 18:2 (2018), 221--238;
\bibitem{B26} A.B. Bogatyrev, “Chebyshev ansatz modes for three band (Stiefel) filter problem”//ArXiv,  2026
\bibitem{Nat}  S.M.Natanzon,  Moduli of surfaces, real algebraic curves and and their superanalogs -- AMS Translation of Math. Monographs, 2004
\bibitem{FK} H.M. Farkas and I.Kra,  Riemann Surfaces
-- Springer Verlag, NY, Heidelberg, Berlin, 1980
\bibitem{DHB} Deconinck,~B.; Heil,~M.; Bobenko,~A.; van Hoeij,~M.; and Schmies,~M. Computing Riemann Theta Functions// Math. Comput. 73, 1417-1442, 2004.
\bibitem{DT} L.N. Trefethen, T.A. Driscoll, Schwarz-Christoffel Mapping --      Cambridge Univ. Press, 2002.
\bibitem{GH} P.Griffiths, J.Harris, 'Principles of Algebraic geometry', Vol.1,  Wiley, NY, 1994 
\bibitem{Mu} D.Muford,  Tata lectures on theta I,II, Birkh\"auser, 1983. 
\end{thebibliography}
\end{document}